\documentclass{article}
\usepackage[]{epsf,epsfig,amsmath,amssymb,amsfonts,latexsym}
\title{Convex Hulls in the Hyperbolic Space}
\newtheorem{theorem}{Theorem}[section]
\newtheorem{lemma}{Lemma}[section]

\newtheorem{remark}{Remark}[section]

\def\qed{\hfill $\vcenter{\hrule height .3mm
\hbox {\vrule width .3mm height 2.1mm \kern 2mm \vrule width .3mm
height 2.1mm} \hrule height .3mm}$ \bigskip}

\def\RR{\mathbb{R}}
\def\H{\mathbb{H}}

\author{Itai Benjamini and Ronen Eldan\thanks{Partially supported by the Israel Science Foundation and by a Farajun Foundation Fellowship}}
\begin{document}
\maketitle
\begin{abstract}
We show that there exists a universal constant $C>0$ such that the convex hull of any $N$ points in 
the hyperbolic space $\mathbb{H}^n$ is of volume smaller than $C N$, and that for any dimension $n$ there exists
a constant $C_n>0$ such that for any set $A \subset \mathbb{H}^n$,
$$
Vol(Conv(A_1)) \leq C_n Vol(A_1)
$$
where $A_1$ is the set of points of hyperbolic distance to $A$ smaller than $1$.
\end{abstract}

\section{Introduction}
The aim of this paper is to prove two fairly basic facts about the volume of convex sets in the hyperbolic space. It is a well known fact that unlike the Euclidean case, the volume of a simplex in the hyperbolic space is bounded from above. The first goal of this note is to show that the volume of any polytope is sublinear with respect to its number of vertices. In dimension 2 the latter fact follows immediately from the former, since any polytope can be triangulated such that the number of triangles is smaller than the number of vertices. This is argument obviously does not work in higher dimensions, since the number of simplices in a triangulation in not linear with respect to the number of vertices. The authors could not find a simple, clean argument which shows this. \\
Our first result may be formulated as follows:
\begin{theorem} \label{linear}
For each dimension $n$ there exists a constant $C_n$ such that the volume of the convex hull of any $N$ points in
$\mathbb{H}^n$ is smaller than $C_n N$.
\end{theorem}

\begin{remark}
An easy corollary of theorem \ref{linear} is the fact that the volume of the fundamental domain of any reflection group in 
$\H^n$ is bounded by a constant which only depends on the dimension and on the number of generators of the group.
\end{remark}
% Compare to euclidean case...
%
Our second result is an application of the first one for estimating the growth of the volume of a set when taking the convex hull. It is easily seen that the volume of a general set can grow in an arbitrary fashion when taking its convex hull. However, it turns out that for a set of bounded curvature this is not the case. For a set $A \subset \mathbb{H}^n$, the $\epsilon$-extension of $A$ is defined as,
$$
A_\epsilon := \bigcup_{x \in A } B_{H}(x, \epsilon)
$$
Where $B_{H}(x,\epsilon)$ is ball of radius $\epsilon$ around $x$ with respect to the hyperbolic metric.
Our second results reads,
\begin{theorem} \label{cvxhull}
For each dimension $n$, and every $\epsilon > 0$, there exists a constant $C(n, \epsilon)$ such that the following holds: For any measurable $A \subset \H^n$, one has
$$
Vol(Conv(A_\epsilon)) < C(n, \epsilon) Vol(A_\epsilon)
$$
\end{theorem}
It is clear that in Euclidean space, even if one adds the assumption that the set is connected, the above is far from true.

\begin{remark}
In dimension $2$, the above theorem has a fairly simple proof in case the set is connected: Let $A \subset \mathbb{H}^2$ be a connected set. It follows from the isoperimetric inequality that the surface area of $A_1$ is comparable to its volume. Moreover, when taking the convex hull of a connected set in in 2 dimensional riemmanian space, the perimeter becomes smaller. It follows that $$Vol_2(conv(A_1)) < C_1 Vol_1(\partial conv(A_1)) < C_1 Vol_1(\partial A_1) < C_2 Vol(A_1)$$
where the first and last inequalities follow from the isoperimetric inequality of $\mathbb{H}^2$.
\end{remark}

A section will be devoted to the proof of each of the theorems. The main idea of the proof of the first theorem is to compare the volume of a convex polytope with that of the union of cone-like objects centered around the vertices of the polytope. One of its main ingredients is a calculation
which roughly shows that every cone centered at a point in infinity and such that the endpoint of every line
in its boundary has a right angle with the geodesic coming from the origin has a bounded volume. This calculation might be of benefit
in understanding the distribution of mass in hyperbolic convex sets. The latter theorem follows from the former rather easily. 

Before we move on to the proofs, let us introduce some notation. We consider the Klein model for the hyperbolic space $\mathbb{H}^n$. For a detailed construction refer to \cite{CFKP}. The Klein model is the Euclidean unit ball in $\RR^n$, denoted by $B_2^n$, equipped with the following metric:
$$
ds^2_{K} = \frac{dx_1^2 + ... + dx_n^2}{1 - x_1^2 - ... - x_n^2} + \frac{(x_1 dx_1 + ... + x_n dx_n)^2}{(1 - x_1^2 - ... - x_n^2)^2} =
$$
$$
\frac{dx^2}{1 - r^2} + \frac{r^2 dr^2}{(1-r^2)^2}
$$
The volume form on the Klein model has the expression,
$$
v_n(r) = \frac{1}{(1-r^2)^{ \frac{n-1}{2} }} \sqrt {\frac{1}{1-r^2} + \frac{r^2}{(1-r^2)^2}} dx = \frac{1}{(1-r^2)^{ \frac{n+1}{2} }} dx
$$
The main advantage of the Klein model over other models, for our purposes, is the fact that geodesics with respect to the hyperbolic metric are
also geodesics with respect to the Euclidean metric, which means that the hyperbolic convex hull of a set is the same as the Euclidean one. A related work of Rivin, \cite{R} has several applications of this fact. \\
For two points $x,y \in B_2^n$ we denote by $x+y$ the standard sum of $x$ and $y$ with respect to the Euclidean linear structure, by 
$|x-y|$ we denote the Euclidean distance between $x$ and $y$, and finally by $d_H(x,y)$ we denote the hyperbolic distance between $x$ and $y$.

\section{The volume of the convex hull of $N$ points is sublinear}
Let $x_1,...,x_N \in \mathbb{H}^n$, define $A = conv(x_1,...,x_N)$.
Our goal in this section is to give an upper bound for 
$Vol(A)$ which depends linearly on $N$. \\ Clearly, we can assume WLOG that $x_1,...,x_N \in S$, and $x_0 = 0$. 
Furthermore, by applying a slight perturbation and using the continuity of the volume, we can assume that the $n-1$-dimensional facets of the polytope $A$ are all simplices. \\ \\
Let's introduce some notation. Define $S$ to be the sphere at infinity, $S = \mathbb{H}^n (\infty) = \partial B_2^n$. For each $x \in S$ let $T_x$ be the unit sphere of the tangent space of $S$ at the point $x$,
which can be identified with $S \cap x^\perp$.

For each $\theta \in T_x$, let
$A_x(\theta) = A \cap span^+ \{x, \theta \}$ (where $span^+$ denotes the positive span) and let $L_x(\theta)$  be the (unique) line which lies on the relative boundary
of $A_x(\theta)$ and passes through $x$ and not through $0$. 
Let $y_x(\theta)$ be the endpoint of $L_x(\theta) \cap B_2$ which is not $x$, 
and let $z_x(\theta)$ be the endpoint of $L_x(\theta) \cap A$ which is not $x$. \\
We make the following definitions,
$$
C_x(\theta) := conv(x, \frac{x + y_x(\theta)}{2}, 0)
$$
$$
\tilde C_x(\theta) := conv(x, \frac{x + z_x(\theta)}{2}, 0)
$$
(the addition is taken with respect to the euclidean structure)
and,
$$
C_x := \bigcup_{\theta \in T_{x}} C_{x}(\theta), ~~ \tilde C_x := \bigcup_{\theta \in T_{x}} \tilde C_{x}(\theta).
$$
Our proof will consist of two main steps. The first one will be to show that there exists a constant $C_n$ such that
\begin{equation} \label{approxbycones}
Vol(A) \leq C_n \sum_{i=1}^N Vol(\tilde C_{x_i}) \leq C_n \sum_{i=1}^N Vol(C_{x_i})
\end{equation}
(the second inequality is obvious by the fact that $\tilde C_i \subseteq C_i$).
The second step will be to show that the volume of $C_i$ is bounded by a 
constant.

We start with proving (\ref{approxbycones}). To this end, let $F$ be an $n-1$ dimensional 
facet of $A$. Assume WLOG that,
$$
F = conv(x_1,...,x_n)
$$
Denote,
$$
D = conv(0, F), ~~ D_i = D \cap \tilde C_{x_i}
$$
It is clear that (\ref{approxbycones}) will follow immediately from the next lemma:
\begin{lemma}
In the above notation, 
$$
Vol(D) \leq 2^n \sum_{i=1}^n Vol(D_i) 
$$
\end{lemma}
\textbf{Proof:}
We can assume $D$ has a nonempty interior, in which case each $y \in D$, can
be uniquely expressed as $y = \sum_{j=1}^n \alpha_j x_j$, $0 \leq \alpha_j \leq 1$. For each $1 \leq i \leq n$, 
define a function $T_i$:
$$
T_i(\sum_{j=1}^n \alpha_j x_j) = \frac{1}{2} \sum_{j=1}^n \alpha_j x_j +  \frac{1}{2} ( \sum_{j=1, j}^n \alpha_j)   x_i
$$
evidently, $T_i(D) = D_i$. \\
Next we note that for each $y \in D$, there exists $1 \leq i \leq n$ such that $|T_i(y)| \geq |y|$. Indeed, this follows from the fact that the vectors $\{T_1(y) - y, ..., T_n(y) - y \}$ positively span an $n-1$ dimensional subspace, and from the convexity of $| \cdot |$. Define,
$$
\ell(y) = \arg \max_{1 \leq i \leq n} |T_i(y)|
$$
and,
$$
T(y) = T_{\ell(y)}(y)
$$
It is easy to see that $T(y)$ is well defined and differentiable for a set whose complement is of measure $0$ in $D$.
We can now calculate,
$$
Vol(D) = \int_D v_n(|x|) dx \leq \int_D v_n(T(x)) dx \leq
$$
$$
\sum_{i=1}^n \int_D v_n(T_i(x)) dx = \sum_{i=1}^n \det(T_i) \int_{C_{x_i}} v_n(|x|) dx = 2^n \sum_{i=1}^n Vol(D_i).
$$
\qed

We move on to the second step, namely, proving the following:
\begin{lemma}
In the above notation, $C_{x_i}$ is bounded by a constant depending only on $n$.
\end{lemma}
\textbf{Proof:} Note that $C_x(\theta)$ is a right triangle. Denote its angle at the origin by $\varphi$. 
One can calculate the volume of $C_i$ by means of polar integration (to be exact, by means of revolution around the axis $[0,x]$):
$$
Vol(C_x) = Z_n \int_{T_x} \int_{C_x(\theta)} d(y,[0,x])^{n-2} v_n(y) d y d \sigma(\theta)
$$
Where $Z_n$ is some constant depending only on the dimension, $\sigma$ is the haar measure on 
$T_x$, and,
$$
d(y,[0,x]) = \min_{0 \leq t \leq 1} |tx -y|.
$$
Let's pick a coordinate system for $span \{x, \theta \}$ in the following way:
define the origin to be $x$, and let $e_1$ be the unit vector in the direction $-x$, and
$e_2$ be the unit vector in the direction of $\theta$. For 
$u,v \in R^+$ denote $(u,v) = (1-u) x + v \theta$. The volume form
$v_n$ becomes,
$$
v_n(u,v) = \frac{1}{(1 - (1-u)^2 - v^2)^{\frac{n+1}{2}}}.
$$
And we get,
$$
\int_{T_x} \int_{C_x(\theta)} d(y,[x,0])^{n-2} d v_n(y) dy = \int_0^1 \int_0^{L(u)} 
\frac{v^{n-2}}{(1 - (1-u)^2 - v^2)^{\frac{n+1}{2}}} dv du 
$$
where 
$$L(u) = 
\begin{cases}
u \cot \varphi &,\; 0 \leq u \leq \sin^2 \varphi  \\
(1-u) \tan \varphi &,\; \sin^2 \varphi \leq u \leq 1.
\end{cases}
$$
We estmate,
$$
\int_0^1 \int_0^{L(u)} 
\frac{v^{n-2}}{(1 - (1-u)^2 - v^2)^{\frac{n+1}{2}}} dv du \leq \int_0^1 \int_0^{L(u)}  \frac{v^{n-2}}{(2u - u^2 - L(u)^2)^{\frac{n+1}{2}}} dv du = 
$$
$$
\int_0^1 \int_0^{L(u)}  \frac{v^{n-2}}{(u + t(u))^{\frac{n+1}{2}}} dv du
$$
where 
$$t(u) = u - u^2 - L(u)^2.$$
Now, $t(0) = t(1) = 0$ and $t(\sin^2 \varphi) = \sin^2 \varphi - \sin^4 \varphi - \sin^4 \varphi \cot^2 \varphi = 0$ which means that $t(u) \geq 0$ for $0 \leq u \leq 1$. So we have,
$$
\int_0^1 \int_0^{L(u)} 
\frac{v^{n-2}}{(1 - (1-u)^2 - v^2)^{\frac{n+1}{2}}} dv du \leq \int_0^1 \int_0^{L(u)} \frac{v^{n-2}}{u^{\frac{n+1}{2}}} dv du =
$$
$$
\frac{1}{n-1} \int_0^1 \frac{L(u)^{n-1}}{u^{\frac{n+1}{2}}} du = \frac{1}{n-1} \int_0^{\sin^2 \varphi} (\cot \varphi)^{n-1} u^{\frac{n-3}{2}} du + \frac{1}{n-1}\int_{\sin^2 \varphi}^{1} (\tan \varphi)^{n-1} (1-u)^{\frac{n-3}{2}} du
$$
We may assume $\varphi < \arctan(\frac {1}{10})$ by starting with a dense-enough net on the sphere. 
This easily implies that the second summand is smaller than $1$. As for the first summand,
$$
\int_0^{\sin^2 \varphi} (\cot \varphi)^{n-1} u^{\frac{n-3}{2}} du = (\cot \varphi)^{n-1} (\sin \varphi)^{n-1} = \cos^{n-1} \varphi.
$$
So finally, we have,
$$
Vol(C_x) \leq C_n' + \cos^{n-1} \varphi + 1 \leq C_n
$$
\qed \\
As explained in the beginning of this section, joining the results of these two lemmata proves theorem ($\ref{linear}$). \\ \\
We finish the section with a few remarks.
\begin{remark}
If we follow the calculations in this section carefully, we may find that that in fact,
$$
\lim_{n \to \infty} C_n = 0
$$
Therefore, there exists a universal constant $C > 0$ such that the volume of a polytope with $N$
vertices of any dimension is smalelr than $C N$. It is interesting to compare this with the asymptotics of the volume of simplices
found in \cite{M}.
\end{remark}
\begin{remark}
It can easily be seen that if for each $N$ we denote by $V_n(N)$ the
 maximal volume of a polytope with $N$ vertices in $\H^n$, then 
$$V_n(N) > c_n N, ~~ \forall N\geq n+1$$ for some constant
$c_n$ depending only on the dimension. This is easily achieved by partitioning $N$ to subsets of size $n+1$ and
constructing disjoint simplices of volume bounded from below. This fact shows us that up to these constants, our result
is, in some sense, sharp.
\end{remark}
\begin{remark}
In Euclidean space, almost all of the volume of any simplex is very far from its vertices. The calculation carried
out in this section suggests that in hyperbolic space this may not be the case. It may be interesting
to find out if the following assertion is true: do there exist constants $c_n \to 0, r_n \to \infty$ (as
$n \to \infty$) such that for any $n$-dimensional simplex such that the distance between each two vertices is at least $r_n$,
at least $0.9$ of its volume is at distance $< c_n r_n$ from one of its vertices? \\
If the latter is true, the calculation carried out in this section would imply that in some sense, unlike the Euclidean case, 
most of the volume any polytope is rather close to its vertices.

\end{remark}

\section{The convex hull of a set whose boundary has bounded curvature}
It is a well known fact that a simplices in $\H^n$ have volume bounded by some
universal constant. This easily implies the following fact:
\begin{lemma} \label{balls}
For each $n \in \mathbb{N}$, there exists a consant $C_n > 0$ such that the convex hull of the union of any
$n+1$ metric balls of radius $1$ in $\H^n$ has a volume smaller than $C_n$.
\end{lemma}
\textbf{Proof:} This immediately follows from the facts that a ball of radius 1 is contained in the convex hull of finitely many points, and that
the volume of any simplex is bounded by a constant. \qed \\
We are now ready to prove our second result. \\
\textbf{Proof of theorem (\ref{cvxhull}):}
Let $A \subset \mathbb{H}^n$.
In view of (\ref{linear}) and lemma (\ref{balls}), it is clearly
enough to show that there exists a constant $C_n$ depending only
on $n$ such that the covering number of $A_\epsilon$ by $\epsilon$-balls is smaller than $C(n,\epsilon) Vol(A_\epsilon)$. \\
Let $N$ be the maximal packing number of $A$, hence, the maximal number of points $x_1,...,x_N \in A$ such that 
\begin{equation} \label{packing}
d_H(x_i, x_j) > \epsilon, ~~ \forall 1 \leq i < j \leq N.
\end{equation}
By the maximality of the packing, we have,
$$
A \subset \bigcup_{i=1}^N B(x_i, \epsilon)
$$
which implies that,
$$
\bigcup_{i=1}^N B(x_i, \frac{\epsilon}{2}) \subset A_\epsilon \subset \bigcup_{i=1}^N B(x_i, 2 \epsilon)
$$
which shows that $N$ is also the covering number of $A$. Moreover, the last relation and (\ref{packing}) imply that,
\begin{equation} \label{vola}
Vol(A) \geq N Vol (B(0, \frac{\epsilon}{2}))
\end{equation}
which is exactly what we need.

\begin{remark}
It may be interesting to find the correct assymptotics for the optimal constants $C(n, \epsilon)$. In view of lemma (\ref{balls})
it is quite clear that the constant provided by our proof will be far from optimal. It is interesting to ask whether the constants $C(n,1)$ are
bounded by some universal constant, in other words, if the maximal ratio between the volume of a set and its convex hull is universally bounded.
\end{remark}

\bigskip {\noindent Department of Mathematics, Weizmann Institute, Rehovot, Israel \\  {\it e-mail address:}
\verb"itai.benjamini@weizmann.ac.il" }

\bigskip {\noindent School of Mathematical Sciences, Tel-Aviv University, Tel-Aviv
69978, Israel \\  {\it e-mail address:}
\verb"roneneldan@gmail.com" }

\end{document}